\newcommand{\calO}{\mathcal{O}}

\newcommand{\calR}{\mathcal{R}}

\newcommand{\bbC}{\mathbb{C}}

\newcommand{\bbP}{\mathbb{P}}
\newcommand{\bbQ}{\mathbb{Q}}

\newcommand{\bbZ}{\mathbb{Z}}

\newcommand\Pic{{\text{Pic}}}

\newcommand\Weil{{\text{Weil}}}
\newcommand{\RHS}{\text{RHS}}
\newcommand{\LHS}{\text{LHS}}

\newcommand{\Sing}{\text{Sing}}

\newcommand{\disc}{\text{disc}}

\documentclass[11pt]{amsart}
\usepackage{amscd, amsmath, amsthm, amssymb}
\newtheorem{theorem}{Theorem}[section]
\newtheorem{lemma}[theorem]{Lemma}
\newtheorem{proposition}[theorem]{Proposition}
\newtheorem{conjecture}[theorem]{Conjecture}
\newtheorem{corollary}[theorem]{Corollary}
\theoremstyle{definition}     
\newtheorem{definition}[theorem]{Definition}

\newtheorem{example}[theorem]{Example}

\theoremstyle{remark}
\newtheorem{remark}[theorem]{Remark}
\numberwithin{equation}{section}

\begin{document}
\title[rationality criterion]
{A Rationality Criterion for projective surfaces - partial solution to Koll\'ar's Conjecture}

\author[J. Keum]{JongHae Keum }
\address{School of Mathematics, Korea Institute for Advanced
Study, Seoul 130-722, Korea } \email{jhkeum@kias.re.kr}
\thanks{Research supported by KOSEF grant R01-2003-000-11634-0}
\dedicatory{Dedicated to Igor Dolgachev on his sixtieth birthday}
\begin{abstract}
Koll\'ar's conjecture states that a complex projective surface $S$ with quotient
singularities and with $H^2(S,\bbQ)\cong \bbQ$ should be rational if its
smooth part $S^0$ is simply connected.

We confirm the conjecture under the additional condition that the exceptional divisor in
a minimal resolution of $S$ has at most 3 components over each
singular point of $S$.
\end{abstract}
\maketitle
\section{Introduction}

In his study of Seifert structures on simply connected rational
 homology spheres, J\'anos Koll\'ar suggested the following conjecture
 (\cite{Ko} Conjecture 42, or
Conjecture 79 its differential geometric equivalent.):

\begin{conjecture}\label{con} Let $S$ be a projective surface with
quotient singularities such that
\begin{itemize}
\item[(1)] $H^2(S, \bbQ)\cong \bbQ$,
\item[(2)] $\pi_1(S^0)=\{1\}$, where $S^0$ is its smooth part.
\end{itemize}
Then $S$ is rational.
\end{conjecture}

In this paper we confirm the conjecture under the additional condition
that the exceptional divisor in
a minimal resolution of $S$ has at most 3 components over each
singular point of $S$.
More precisely, we prove the following:

\begin{theorem}\label{main} Let $S$ be a projective surface with
quotient singularities such that
\begin{itemize}
\item[(1)] $H^2(S, \bbQ)\cong \bbQ$,
\item[(2)] $H_1(S^0, \bbZ)=0$,
\item[(3)] the inverse image  $f^{-1}(p)$ has at most 3 components for each
singular point $p$ in $S$, where $f: S'\to S$ is a minimal
resolution.
\end{itemize}
Then $S$ is rational.
\end{theorem}

Note that the condition (2) $H_1(S^0, \bbZ)=0$ is weaker than
$\pi_1(S^0)=\{1\}$.

We also remark that if $S$ is non-singular and satisfies the
conditions (1) and (2) of Theorem  \ref{main}, then $S$ is
 either the complex projective plane $\bbC\bbP^2$
   or a surface of general type
with $q=p_g=0$,  $3c_2=c_1^2=9$, so called a fake projective
plane. Recently G. Prasad and S.-K. Yeung have shown that no fake
projective plane with $H_1(S, \bbZ)=0$ exists \cite{PY}.

Throughout this paper, we work over the field $\bbC$ of complex numbers.

\medskip
{\it Acknowledgements.} I like to thank J\'anos Koll\'ar for useful
conversations through e-mails. I am also grateful to the referee for
many helpful comments.

\section{Preliminaries}

\begin{lemma}\label{1} Let $V$ be an irreducible reduced complex analytic
  space, and $V^0$ its smooth part. Let $f: V'\to V$ be a
resolution of singularities.  Then
\begin{itemize}
\item[(1)] The inclusion $V^0\subset V'$ gives surjective homomorphisms
$$\pi_1(V^0)\to \pi_1(V'), \quad H_1(V^0, \bbZ)\to H_1(V', \bbZ).$$
\item[(2)] If $f$ has connected fibres, it induces surjective homomorphisms
$$\pi_1(V')\to \pi_1(V), \quad H_1(V', \bbZ)\to H_1(V, \bbZ).$$
\end{itemize}
\end{lemma}

\begin{proof} It suffices to prove the assertions for fundamental groups.

The first assertion follows from the fact that
the complement $V'\setminus V^0$ has real codimension $\ge 2$.

If $f$ has connected fibres, a loop in $V$ can be lifted to a loop in
$V'$.
\end{proof}

When the resolution is projective, the condition in (2) is always
satisfied by Zariski's Main Theorem.

The following also can be proved by a standard argument.
For a proof, we refer the reader, e.g. to \cite{Ko} Proposition 40.

\begin{proposition}\label{2} Let $S$ be a projective surface with
quotient singularities such that $H^2(S, \bbQ)\cong \bbQ$ and  $H^1(S,
 \bbQ)=0$. Let $f: S'\to S$ be a
resolution of singularities. Then
\begin{itemize}
\item[(1)] $H^1(S, \calO_S)=H^1(S', \calO_{S'})=0$.
\item[(2)] $H^2(S, \calO_S)=H^2(S', \calO_{S'})=0$.
\end{itemize}
\end{proposition}

Recall the definition of {\it Kodaira (logarithmic) dimension}.
Let $V^0$ be a nonsingular variety and let
$V$ be a {\it smooth completion} of $V^0$, i.e., $V$ is nonsingular projective and
$D := V \setminus V^0$ is an integral reduced divisor with simple normal crossings.
If $H^0(V, m(K_V+D)) = 0$ for all $m \ge 1$, the
{\it Kodaira (logarithmic) dimension} $\kappa(V^0) = -\infty$.
Otherwise, $|m(K_V+D)|$ gives rise to a rational map $\varphi_m$
for some $m$ and the {\it Kodaira dimension}
$\kappa(V^0)$ is the maximum of $\dim(\varphi_m(V^0))$.

\par
The Kodaira dimension of $V^0$ does not depend on the choice of
the completion $V$ \cite{Iitaka}. Also
$\kappa(V^0)$ takes value in $\{-\infty, 0, 1, \dots, \dim V^0\}$.

\par
Obviously, $\kappa(V)\le \kappa(V^0)$.

\begin{proposition}\label{3} Let $S$ be a projective surface with
quotient singularities such that $H^2(S, \bbQ)\cong \bbQ$ and
$H_1(S^0,\bbZ)=0$. Let $f: S'\to S$ be a minimal
resolution. Then one of the following cases occurs.
\begin{itemize}
\item[(1)] $S$ is rational.
\item[(2)] $S'$ is a surface, not necessarily minimal,
with $q=p_g=0$, $\kappa(S')=1$, $H_1(S',\bbZ)=0$ and $\kappa(S^0)=2$.
\item[(3)] $S'$ is a surface of general type,
not necessarily minimal, with $q=p_g=0$,  $H_1(S',\bbZ)=0$.
\end{itemize}
\end{proposition}

\begin{proof} Since $H_1(S^0,\bbZ)=0$,
$H_1(S',\bbZ)=H_1(S,\bbZ)=0$ by Lemma \ref{1}. By Proposition \ref{2}, $S'$ is a
surface with $q=p_g=0$, $H_1(S',\bbZ)=0$. By classification theory (see \cite{BHPV}),
either $S'$ is rational or $\kappa(S')\ge 1$. Since $H^2(S, \bbQ)\cong
\bbQ$, $S$ has Picard
number 1 and $H^2(S,
\bbQ)$ is positive definite. Thus for every curve $C$ on $S$,  $C^2\ge 0$.
In particular, $S$ is relatively minimal, i.e.
there is no curve $C$ with $K_S\cdot C<0$, $C^2<0$.

Assume $\kappa(S')=1$.
If $\kappa(S^0)=1$, then there is an elliptic fibration on $S$
(\cite{Ka1} Theorem 2.3 or \cite{Mi} Ch.II Theorem 6.1.4 or \cite{KZ}
Theorem 4.1.), thus $\Pic(S)$ has rank at least 2,
a contradiction.
Thus $\kappa(S^0)=2$ and the assertion follows.
\end{proof}

\medskip

Replacing the condition $H_1(S^0,\bbZ)=0$ by $\pi_1(S^0)=\{1\}$,
one gets the following

\begin{corollary}\label{4} Let $S$ be a projective surface with
quotient singularities such that $H^2(S, \bbQ)\cong \bbQ$ and
$\pi_1(S^0)=\{1\}$. Let $f: S'\to S$ be a minimal
resolution. Then one of the following cases occurs.
\begin{itemize}
\item[(1)] $S$ is rational.
\item[(2)] $S'$ is a simply connected surface, not necessarily minimal,
with $q=p_g=0$, $\kappa(S')=1$, and $\kappa(S^0)=2$.
\item[(3)] $S'$ is a simply connected surface of general type,
not necessarily minimal, with $q=p_g=0$.
\end{itemize}
\end{corollary}

So far, no example $S$ satisfying the condition of Corollary \ref{4}
and belonging to the cases (2) or (3) has been
found, and it is not likely such an example exists. On this basis
J\'anos Koll\'ar suggests his conjecture (\cite{Ko} 41).

\begin{lemma}\label{5} Let $S$ be a normal compact surface with
  rational singularities, and $f: S'\to S$ a
resolution of singularities. Let $\calR\subset H^2(S', \bbZ)$ be the
subgroup generated by the cohomology classes of the exceptional
curves of $f$. Set
 $$\overline{\calR}:=\{\alpha\in H^2(S',\bbZ): m\alpha\in \calR \quad
 {\rm for} \quad {\rm some} \quad m\in \bbZ\setminus \{0\} \}.$$
Then the following are equivalent
\begin{itemize}
\item[(1)] $H_1(S^0,\bbZ)=0$.
\item[(2)] $q(S')=0$ and $\calR=\overline{\calR}$.
\end{itemize}
\end{lemma}

\begin{proof}
Assume (1). By Lemma \ref{1}, $H_1(S', \bbZ)=0$, and hence
$q(S')=0$. By the universal coefficient
theorem, $H^2(S', \bbZ)$ is torsion free, so is $\overline{\calR}$.
Since  $q(S')=0$, $\Pic(S')$ can be regarded as a primitive sublattice
of  $H^2(S', \bbZ)$. Since $\calR\subset \Pic(S')$,
$\overline{\calR}\subset \Pic(S')$. If $\calR\neq \overline{\calR}$,
then there would be a finite \'etale cover of $S^0$, thus
$H_1(S^0,\bbZ)\neq 0$, a contradiction.

\medskip
Assume (2). Since $q(S')=0$, $\Pic(S')$
embeds in $H^2(S', \bbZ)$.
Suppose $H_1(S^0,\bbZ)\neq 0$. Then there is a finite \'etale cover of
$S^0$, thus there exist an element $\alpha\in \Pic(S')$ and an
integer $m>1$ such that $m\alpha$ is either trivial or linearly equivalent to an
effective divisor supported in the exceptional
set of $f$, but $\alpha$ is not. Since $\Pic(S')\subset H^2(S',
\bbZ)$, this implies that $\overline{\calR}\neq\calR$.
\end{proof}

\begin{definition} Let $p\in F$ be a normal surface singularity. Then $F$
 is a cone over a real
  3-manifold $M$ called the {\it link}.

If the singularity is
 rational, then $H^1(M,\bbZ)=0$ and $H^2(M,\bbZ)$ is torsion.
\end{definition}

For surfaces with $H^2(S,\bbQ)=\bbQ$, J. Koll\'ar gives more precise
information in terms of links.

\begin{proposition}\label{6} $($\cite{Ko} Corollary $43)$ Let $S$ be a normal compact surface with
  rational singularities $p_i$ with links $M_i$. Assume that
$H_1(S,\bbZ)=0$ and  $H^2(S,\bbQ)=\bbQ$. Then the following are
equivalent
\begin{itemize}
\item[(1)] $H_1(S^0,\bbZ)=0$.
\item[(2)] The Weil divisor class group $\Weil (S)\cong \bbZ$.
\item[(3)] Each $H^2(M_i,\bbZ)$ is cyclic, their orders $m_i$ are
  pairwise coprime and there is a Weil divisor $B'$ which generates
 $H^2(M_i,\bbZ)$ for every $i$.
\item[(4)] There is a Weil divisor $B$ with $B^2=1/\Pi m_i$.
\item[(5)] There is a Cartier divisor $H$ and a Weil divisor $B$ with
  $H^2=\Pi m_i$, $B\cdot H=1$.
\end{itemize}
\end{proposition}

The folowing result due to Y. Miyaoka plays a crucial role in the
proof of our main theorem.

\begin{theorem}\label{m} $($\cite{Miya} Theorem $1.1)$ Let $S$ be a projective surface with
quotient singularities. Denote  by $\Sing(S)$ the set of singular points
of $S$.  Let $f: S'\to S$ be a minimal
resolution and $E$ be the inverse image $f^{-1}(\Sing(S))$, a
reduced integral divisor. Assume $K_{S'}+E$ has Zariski decomposition with
positive part $P$ and negative part $N+N'$, where $N$ is supported away
from $E$ and $N'$ is supported in $E$. Then we have the inequality
$$\sum_{p\in \Sing(S)} (e(E_p)-\frac{1}{|G_p|})\le
c_2(S')-\frac{1}{3}P^2-\frac{1}{4}N^2,$$
where $e(E_p)$ is the Euler number of  $E_p:=f^{-1}(p)$ and $G_p$ is
the local fundamental group of $p$.
\end{theorem}

\begin{corollary}\label{corM} Let $S$ be  a projective surface with
quotient singularities, and $f: S'\to S$ be a minimal
resolution. Assume $S$ is relatively minimal, i.e. there
 is no curve $C$ with $K_S\cdot C<0$, $C^2<0$. Assume  $\kappa(S^0)\ge
 0$. Then we have the inequality
\begin{equation}\label{ineq1}
\sum_{p\in \Sing(S)} (e(E_p)-\frac{1}{|G_p|})\le
c_2(S')-\frac{1}{3}K_S^2.
\end{equation}
\end{corollary}

\begin{proof} The canonical divisor $K_S$ is numerically effective by
\cite{MT}, Theorem 2.11 or \cite{KZ}, Theorem 2.1.
Since a quotient singularity is just a log terminal singularity, we
have
$$K_{S'}+E=f^*K_S+(E-D),$$
where $E-D$ is an effective $\bbQ$-divisor whose support is equal to $E$.
Thus the positive part $P$ of Zariski decomposition of $K_{S'}+E$ is $f^*K_S$ and
the negative part is $E-D$.
\end{proof}

Let $S$ be  a projective surface with quotient singularities.
Then one can write
\begin{equation}\label{K}
K_{S'}=f^*K_S-\sum_{p\in \Sing(S)} D_p
\end{equation}
where $D_p$ is an effective  $\bbQ$-divisor supported in $E_p=f^{-1}(p)$.

\begin{corollary}\label{corM2} Let $S$ be a projective surface
with quotient singularities such that $H^2(S, \bbQ)\cong \bbQ$
and $H_1(S^0, \bbZ)=0$. Let $f: S'\to S$ be a minimal
resolution.
Assume that $S$ is not rational.
 Then we have the inequality
\begin{equation}\label{ineq}
\sum_{p\in \Sing(S)} e(E_p)-c_2(S')+\frac{1}{3}K_{S'}^2\le
\sum_{p\in \Sing(S)}(\frac{1}{|G_p|}+\frac{1}{3}D_p^2).
\end{equation}
\end{corollary}

\begin{proof} By Proposition \ref{3},  $\kappa(S^0)=2$. Since
$H^2(S, \bbQ)\cong \bbQ$, $S$ has Picard number 1 and $H^2(S,
\bbQ)$ is positive definite. Thus for every curve $C$ on $S$,  $C^2\ge 0$.
In particular, $S$ is relatively minimal. (This also follows from
Kawamata's Cone Theorem. Indeed, the existence of
a curve $C$ with $K_S\cdot C<0$, $C^2<0$ would imply the existence of
an extremal contraction, which is either divisorial or gives a
fibration, both contradicting to the fact that  $S$ has Picard number
1.) By Corollary \ref{corM},
we get the inequality \eqref{ineq1}.
It remains to see that
$$K_S^2=(f^*K_S)^2=K_{S'}^2-\sum_{p\in \Sing(S)} D_p^2.$$
\end{proof}

\begin{corollary}\label{corM3} Let $S$ be a projective surface
with quotient singularities such that $H^2(S, \bbQ)\cong \bbQ$
and $H_1(S^0, \bbZ)=0$.
Assume that $S$ is not rational. Then $S$ has at most $4$ singular points.
\end{corollary}

\begin{proof}  Let $f: S'\to S$ be a minimal
resolution. By Corollary \ref{corM},
we have the inequality \eqref{ineq1}.
Note that
$$\sum_{p\in \Sing(S)} e(E_p)=b_2(S')-1+|\Sing(S)|=c_2(S')-3+|\Sing(S)|.
$$
Thus the inequality \eqref{ineq1} becomes
\begin{equation}\label{ineq4}
|\Sing(S)|-3\le \sum_{p\in \Sing(S)}\frac{1}{|G_p|}-\frac{1}{3}K_S^2.
\end{equation}
Let $p_1,..., p_r$ be the singular points of $S$, and let $M_1, ..., M_r$ be the
corresponding links. Since the singularities are rational,
$H^2(M_i,\bbZ)$ is isomorphic to the
abelianization of the local fundamental group $G_{p_i}$.
By Proposition \ref{6},
$H^2(M_i,\bbZ)$ is cyclic, and their orders $m_i$ are
  pairwise coprime. Let us assume that $m_1<m_2< ...<m_r$.

Assume $r=|\Sing(S)|\ge 5$. If $m_1>1$, then $m_i$ is greater than
equal to the $i$-th prime number, thus
\begin{equation}\label{ineq5}
\sum_{i=1}^r\frac{1}{|G_{p_i}|}\le \sum_{i=1}^r\frac{1}{m_i}\le
r-3.
\end{equation}
If $m_1=1$, then $\frac{1}{|G_{p_1}|}< \frac{1}{2}$, thus
\begin{equation}\label{ineq6}
\sum_{i=1}^r\frac{1}{|G_{p_i}|}\le \frac{1}{|G_{p_1}|}+ \sum_{i=2}^r\frac{1}{m_i}\le
\frac{1}{2}+(r-\frac{7}{2})=r-3.
\end{equation}
Since $K_S^2>0$, both \eqref{ineq5} and \eqref{ineq6} lead to a contradiction to the inequality
\eqref{ineq4}.
Thus $|\Sing(S)|\le 4$.
\end{proof}

\begin{remark}\label{r1} In the situation of Corollary \ref{corM3}, if
 $|\Sing(S)|=4$, then two of the four singularities have the local fundamental
 group of order 2 and 3, respectively.
\end{remark}

\section{Proof of Main Theorem}

In this section, we prove Theorem \ref{main}.

Fix a singular point $p\in S$, and let $E_{1}, ...,
E_{k}$ ($k\le 3$) be the irreducible components of $E_p=f^{-1}(p)$.
They form a string of smooth rational curves
$$(-n_1) \quad {\rm or}\quad (-n_1)-(-n_2)  \quad {\rm or}\quad (-n_1)-(-n_2)-(-n_3)$$
where $E_{j}$ is a $(-n_j)$-curve. Write
$$D_p=\sum_{j=1}^k a_{j}E_{j}.$$
Note that $0\le a_j<1$.

To use the inequality \eqref{ineq}, we need to estimate $\frac{1}{|G_p|}+\frac{1}{3}D_p^2$.

\begin{lemma}\label{2.1} Fix $p\in \Sing(S)$. Assume that  $f^{-1}(p)$
has $3$ components $E_1$, $E_2$, $E_3$ with $E_i^2=-n_i$.
Assume that $n_1+n_2+n_3\ge 10$. Then
$$\frac{1}{|G_p|}+\frac{1}{3}D_p^2< -\frac{1}{2}.$$
\end{lemma}

\begin{proof} Since $E_j$ is a $(-n_j)$-curve, $E_j\cdot
  K_{S'}=n_j-2$. Intersecting $E_j$ with $f^*K_S$ from \eqref{K}, we see that
$$n_1-2=a_1n_1-a_2$$
$$n_2-2=-a_1+a_2n_2-a_3$$
$$n_3-2=-a_2+a_3n_3$$
Adding the equations, we get
$$\sum n_j-6=\sum a_jn_j-\sum a_j-a_2,$$
hence
$$\sum a_j(n_j-2)=\sum n_j-6-a_1-a_3>\sum n_j-8 \ge 2.$$
Since $D_p^2=-D_p\cdot K_{S'}=-\sum a_j(n_j-2)$, we have
$$\frac{1}{|G_p|}+\frac{1}{3}D_p^2=\frac{1}{|G_p|}-\frac{1}{3}\sum
a_j(n_j-2)<\frac{1}{|G_p|}-\frac{2}{3}\le -\frac{1}{2}.$$
\end{proof}

For the cases where $n_1+n_2+n_3\le 9$, we give an exact estimate
in Table 1.

\medskip

\begin{table}[h]
\begin{center}
\begin{tabular}{|l||r|r|r|r|}
\hline
$(n_1, n_2, n_3)$&$|G_p|$&$(a_1, a_2, a_3)$&$D_p^2$&$\frac{1}{|G_p|}+\frac{1}{3}D_p^2$\\ \hline
$(2,2,2)$&$4$&$(0,0,0)/4$&$0$&$1/4$\\ \hline
$(2,2,3)$&$7$&$(1,2,3)/7$&$-3/7$&$0$\\ \hline
$(2,3,2)$&$8$&$(2,4,2)/8$&$-4/8$&$-1/24$\\ \hline
$(2,2,4)$&$10$&$(2,4,6)/10$&$-12/10$&$-3/10$\\ \hline
$(2,4,2)$&$12$&$(4,8,4)/12$&$-16/12$&$-13/36$\\ \hline
$(3,2,3)$&$12$&$(6,6,6)/12$&$-12/12$&$-1/4$\\ \hline
$(2,3,3)$&$13$&$(4,8,7)/13$&$-15/13$&$-4/13$\\ \hline
$(2,2,5)$&$13$&$(3,6,9)/13$&$-27/13$&$-8/13$\\ \hline
$(2,5,2)$&$16$&$(6,12,6)/16$&$-36/16$&$-11/16$\\ \hline
$(3,2,4)$&$17$&$(9,10,11)/17$&$-31/17$&$-28/51$\\ \hline
$(2,3,4)$&$18$&$(6,12,12)/18$&$-36/18$&$-11/18$\\ \hline
$(2,4,3)$&$19$&$(7,14,11)/19$&$-39/19$&$-12/19$\\ \hline
$(3,3,3)$&$21$&$(12,15,12)/21$&$-39/21$&$-4/7$\\ \hline

\end{tabular}
\end{center}
\caption{}
\end{table}

\begin{lemma}\label{2.2} Fix $p\in \Sing(S)$. Assume that  $f^{-1}(p)$
has $2$ components $E_1$, $E_2$ with $E_i^2=-n_i$.
Assume that $n_1+n_2\ge 8$. Then
$$\frac{1}{|G_p|}+\frac{1}{3}D_p^2< -\frac{1}{2}.$$
\end{lemma}

\begin{proof}  Intersecting $E_j$ with $f^*K_S$ from \eqref{K}, we see that
$$n_1-2=a_1n_1-a_2$$
$$n_2-2=-a_1+a_2n_2$$
Adding the equations, we get
$$\sum n_j-4=\sum a_jn_j-\sum a_j,$$
hence
$$\sum a_j(n_j-2)=\sum n_j-4-a_1-a_2>\sum n_j-6 \ge 2.$$
Thus, we have
$$\frac{1}{|G_p|}+\frac{1}{3}D_p^2=\frac{1}{|G_p|}-\frac{1}{3}\sum
a_j(n_j-2)<\frac{1}{|G_p|}-\frac{2}{3}\le -\frac{1}{2}.$$
\end{proof}

For the cases where $n_1+n_2\le 7$, we give an exact estimate
in Table 2.

\medskip

\begin{table}[h]
\begin{center}
\begin{tabular}{|l||r|r|r|r|}
\hline
$(n_1, n_2)$&$|G_p|$&$(a_1, a_2)$&$D_p^2$&$\frac{1}{|G_p|}+\frac{1}{3}D_p^2$\\ \hline
$(2,2)$&$3$&$(0,0)/3$&$0$&$1/3$\\ \hline
$(2,3)$&$5$&$(1,2)/5$&$-2/5$&$1/15$\\ \hline
$(2,4)$&$7$&$(2,4)/7$&$-8/7$&$-5/21$\\ \hline
$(3,3)$&$8$&$(4,4)/8$&$-8/8$&$-5/24$\\ \hline
$(2,5)$&$9$&$(3,6)/9$&$-18/9$&$-5/9$\\ \hline
$(3,4)$&$11$&$(6,7)/11$&$-20/11$&$-17/33$\\ \hline

\end{tabular}
\end{center}
\caption{}
\end{table}

\begin{lemma}\label{2.3}  Fix $p\in \Sing(S)$. Assume that  $f^{-1}(p)$
has $1$ component $E_1$ with $E_1^2=-n$.
Let $d=|G_p|$. Then $d=n$ and
$$\frac{1}{|G_p|}+\frac{1}{3}D_p^2=\frac{1}{d}-\frac{(d-2)^2}{3d}$$
which equals to $\frac{1}{2}$ if
$d=2$, to $\frac{2}{9}$ if $d=3$, and  $\le -\frac{1}{12}$ if $d\ge 4$.
\end{lemma}

\begin{proof}
$D_p=\frac{d-2}{d}E_1$.
\end{proof}

If $H_1(S^0,\bbZ)=0$, then by Lemma \ref{1}, $H_1(S',
\bbZ)=0$, hence $H^2(S', \bbZ)$ is torsion free and becomes a lattice
with intersection pairing.

From Lemma \ref{5} and Proposition \ref{6}, we also have the following:

\begin{lemma}\label{2.4} Let $S$ be a projective surface with
quotient singularities
satisfying the conditions $(1)$ and $(2)$ of Theorem \ref{main}.
 Write $\calR=\oplus_p \calR_p$ where $\calR_p$ is the sublattice of
 $H^2(S', \bbZ)$ generated by the components of $E_p=f^{-1}(p)$. Then
\begin{itemize}
\item[(1)] The numbers $|G_p/[G_p,G_p]|=|\det\calR_p|$ are pairwise coprime.
\item[(2)] There is an integer $m$ such that $|\det\calR|(f^*K_S)^2=m^2$.
\end{itemize}
\end{lemma}

\begin{proof} Here we give a short proof.

By Lemma \ref{5}, $\calR$ is a primitive sublattice of $H^2(S', \bbZ)$.
Since  $H^2(S', \bbZ)$ is unimodular, we have an isomorphism between
the discriminant groups
$$\disc\calR=\oplus_p \disc\calR_p \cong -\disc\calR^{\perp}.$$
Since $\calR^{\perp}$ is of rank 1, $\disc\calR$ must be cyclic. This
proves (1).

The divisor  $(\det\calR)f^*K_S$ is an integral divisor belonging to
$\calR^{\perp}$, hence $(\det\calR)f^*K_S=mv$ for some integer $m$,
where $v$ is a generator of $\calR^{\perp}$.
Since $v^2=|\det\calR|$, (2) follows.
\end{proof}

From now on,  $S$ denotes a projective surface
satisfying the condition of Theorem \ref{main},
i.e. $S$ is a singular projective surface with
quotient singularities such that
\begin{itemize}
\item[(1)] $H^2(S, \bbQ)\cong \bbQ$,
\item[(2)] $H_1(S^0, \bbZ)=0$,
\item[(3)] the inverse image  $f^{-1}(p)$ has at most 3 components for each
singular point $p$ in $S$, where $f: S'\to S$ is a minimal
resolution.
\end{itemize}
To get a contradiction, we also assume
\begin{itemize}
\item[(4)] $S$ is not rational.
\end{itemize}

\bigskip

In this situation, by Corollary \ref{corM2}, we have the inequality
\eqref{ineq}. By the assumption (3), all singularities of $S$
are cyclic.

We denote the left hand side and the right hand side of the inequality
\eqref{ineq} by
$$\LHS:= \sum_{p\in \Sing(S)} e(E_p)-c_2(S')+\frac{1}{3}K_{S'}^2$$
$$\RHS:= \sum_{p\in \Sing(S)}(\frac{1}{|G_p|}+\frac{1}{3}D_p^2).$$

\begin{lemma}\label{re} Let $S$ be a projective surface with
quotient singularities satisfying the conditions $(1)-(4)$.
Assume that the number of singular points $|\Sing (S)|\ge 3$.
Then $\RHS\le \frac{9}{10}$.
\end{lemma}

\begin{proof} From Lemma \ref{2.1}-\ref{2.3}, we see that
$\frac{1}{|G_p|}+\frac{1}{3}D_p^2 > 0$
for only one of the five types of singularities
$(2,2,2), (2,2), (2,3), (2), (3)$.
Also, by Lemma \ref{2.4}, the pair $(2,2,2)$ and $(2)$ do not occur
simultaneously. Neither the pair $(2,2)$ and $(3)$.

If $|\Sing (S)|\ge 3$, $\RHS$ takes its maximum value
when  $\Sing (S)=(2)+(2,2)+(2,3)$, hence
$$\RHS\le \frac{1}{2}+\frac{1}{3}+\frac{1}{15}=\frac{9}{10}.$$
\end{proof}

{\bf Proof of Theorem \ref{main}.} To get a contradiction, assume
that $S$ is not rational. By Proposition \ref{3} it suffices to
rule out the two cases
\begin{itemize}
\item[(2)] $S'$ is a surface, not necessarily minimal,
with $q=p_g=0$, $\kappa(S')=1$, $H_1(S',\bbZ)=0$ and $\kappa(S^0)=2$.
\item[(3)] $S'$ is a surface of general type,
not necessarily minimal, with $q=p_g=0$,  $H_1(S',\bbZ)=0$.
\end{itemize}

\medskip

Since $q=p_g=0$, by Bogomolov-Miyaoka-Yau inequality (or by Theorem
\ref{m}) we have $c_2(S')\ge 3$.

If $c_2(S')=3$, then $S'=S$ and $S$ is non-singular, hence $S$ is
either the complex projective plane
   or a surface of general type
with $q=p_g=0$,  $3c_2=c_1^2=9$, so called a fake projective
plane. The latter surface has the unit ball in $\bbC^2$ as its
universal covering (this follows from the solution of S.-T. Yau
\cite{Yau} to Calabi conjecture) hence has an infinite fundamental
group. G. Prasad and S.-K. Yeung \cite{PY} have shown that no fake
projective plane with $H_1(S, \bbZ)=0$ exists.

Thus we may assume that $c_2(S')\ge 4$ and $S$ is singular.

By Corollary \ref{corM2}, we have the inequality
\eqref{ineq}. By Corollary \ref{corM3}, we also have $|\Sing(S)|\le
4$.

\medskip\noindent
Case 1. $c_2(S')=4$ and $K_{S'}^2=8$.

\medskip\noindent
In this case $|\Sing(S)|=1$ and $\LHS=2-4+\frac{8}{3}=\frac{2}{3}$, while
$\RHS\le \frac{1}{2}$, a contradiction.

\medskip\noindent
Case 2. $c_2(S')=5$ and $K_{S'}^2=7$.

\medskip\noindent
If $|\Sing(S)|=2$, then $\LHS=4-5+\frac{7}{3}=\frac{4}{3}$, while
$\RHS\le \frac{1}{2}+\frac{2}{9}$, a contradiction.

\noindent
If $|\Sing(S)|=1$, then $\LHS=3-5+\frac{7}{3}=\frac{1}{3}$, while
$\RHS\le \frac{1}{3}$, with equality only when
$(n_1, n_2)=(2, 2)$. In this case $\det\calR=3$ and
$(f^*K_S)^2=K_{S'}^2=7$, a contradiction to Lemma \ref{2.4}(2).

\medskip\noindent
Case 3. $c_2(S')=6$ and $K_{S'}^2=6$.

\medskip\noindent
If $|\Sing(S)|=3$, then $\LHS=6-6+\frac{6}{3}=2$, contradicts to Lemma \ref{re}.

\noindent
If $|\Sing(S)|=2$, then $\LHS=5-6+\frac{6}{3}=1$, while
 $\RHS\le \frac{1}{2}+\frac{1}{3}$, a contradiction.

\noindent
If $|\Sing(S)|=1$, then $\LHS=4-6+\frac{6}{3}=0$, hence we must have
$(n_1, n_2, n_3)=(2, 2, 2)$ or $(2,2,3)$. In the first case (resp. the
second) $\det\calR=4$ (resp. $7$) and $(f^*K_S)^2=6$
(resp. $\frac{45}{7}$).
Both  contradict to Lemma \ref{2.4}(2).

\medskip\noindent
Case 4. $c_2(S')=7$ and $K_{S'}^2=5$.

\medskip\noindent
If $|\Sing(S)|\ge 3$, then $\LHS\ge 7-7+\frac{5}{3}=\frac{5}{3}$,
contradicts to Lemma \ref{re}.

\noindent
If $|\Sing(S)|=2$, then $\LHS=6-7+\frac{5}{3}=\frac{2}{3}$, hence
$(n_1, n_2, n_3)+(n_4)=(2, 2, 2)+(2)$ (no possible combination of type
$(n_1, n_2)+(n_3,n_4)$). In this case $\det\calR_{p_1}=4$ and
$\det\calR_{p_2}=2$, contradicting to Lemma \ref{2.4}(1).

\medskip\noindent
Case 5. $c_2(S')=8$ and $K_{S'}^2=4$.

\medskip\noindent
If $|\Sing(S)|\ge 3$, then $\LHS\ge 8-8+\frac{4}{3}=\frac{4}{3}$,
contradicts to Lemma \ref{re}.

\noindent
If $|\Sing(S)|=2$, then $\LHS=7-8+\frac{4}{3}=\frac{1}{3}$, hence
$(n_1, n_2, n_3)+(n_4, n_5)=(2, 2, 2)+(2,2)$ or $(2, 2, 3)+(2,2)$.
In the first case (resp. the
second) $\det\calR=12$ (resp. $21$) and $(f^*K_S)^2=4$
(resp. $\frac{31}{7}$).
Both  contradict to Lemma \ref{2.4}(2).

\medskip\noindent
Case 6. $c_2(S')=9$ and $K_{S'}^2=3$.

\medskip\noindent
If $|\Sing(S)|\ge 3$, then  $\LHS\ge 9-9+\frac{3}{3}=1$,
contradicts to Lemma \ref{re}.

\noindent
If $|\Sing(S)|=2$, then $\LHS=8-9+\frac{3}{3}=0$, hence
$(n_1, n_2, n_3)+(n_4, n_5,n_6)=(2, 2, 2)+(2,2,3)$ or $(2, 2,
2)+(2,3,2)$ or $(2, 2, 2)+(3,2,3)$.
In the first case, $\det\calR=28$ and $(f^*K_S)^2=\frac{24}{7}$,
contradicts to Lemma \ref{2.4}(2).
In the second case, $\det\calR=4\cdot 8$, and in the third,
$\det\calR=4\cdot 12$, both contradict to Lemma \ref{2.4}(1).

\medskip\noindent
Case 7. $c_2(S')=10$ and $K_{S'}^2=2$.

\medskip\noindent
If $|\Sing(S)|=4$, then $\LHS=
11-10+\frac{2}{3}=\frac{5}{3}$,
contradicts to Lemma \ref{re}.

\noindent
If $|\Sing(S)|=3$, then $\LHS=\frac{2}{3}$, hence
$(n_1, n_2, n_3)+(n_4, n_5,n_6)+(n_7)=(2, 2, 2)+(2,2,3)+(2)$ or
$(2,2,2)+(2, 3,2)+(2)$
(no possible combination of type $(n_1, n_2, n_3)+(n_4, n_5)+(n_6, n_7)$).
In the first case, $\det\calR=56$ and $(f^*K_S)^2=K_{S'}^2-\sum
D_p^2=2+\frac{3}{7}=\frac{17}{7}$,
contradicts to Lemma \ref{2.4}(2).
In the second case, $\det\calR=4\cdot 8\cdot 2$, contradicts to Lemma \ref{2.4}(1).

\medskip\noindent
Case 8. $c_2(S')=11$ and $K_{S'}^2=1$.

\medskip\noindent
If $|\Sing(S)|=4$, then $\LHS= 12-11+\frac{1}{3}=\frac{4}{3}$,
contradicts to Lemma \ref{re}.

\noindent
If $|\Sing(S)|=3$, then $\LHS=\frac{1}{3}$, hence by Lemma \ref{2.4}(1),
$(n_1, n_2, n_3)+(n_4, n_5,n_6)+(n_7,n_8)=(2, 2, 2)+(2,2,3)+(2,2)$.
Then, $\det\calR=4\cdot 7\cdot 3$ and $(f^*K_S)^2=K_{S'}^2-\sum
D_p^2=1+\frac{3}{7}=\frac{10}{7}$,
contradicts to Lemma \ref{2.4}(2).

\medskip\noindent
Case 9. $c_2(S')=3s+3$ and $K_{S'}^2=9-3s$ ($s\ge 3$).

\medskip\noindent
In this case  $|\Sing(S)|\ge (b_2(S')-1)/3=3s/3=s$.

Assume $|\Sing(S)|\ge s+1$, then
\begin{eqnarray*}
\LHS&=& (b_2(S')-1)+|\Sing(S)|-c_2(S')+\frac{1}{3}K_{S'}^2\\
&\ge& 3s+(s+1)-(3s+3)+\frac{1}{3}(9-3s)= 1.
\end{eqnarray*}
Since $s\ge 3$, this contradicts to Lemma \ref{re}.

Assume $|\Sing(S)|=s$, then
$$\LHS= (b_2(S')-1)+|\Sing(S)|-c_2(S')+\frac{1}{3}K_{S'}^2= 0.$$
Since $b_2(S')-1=3s$, $\Sing(S)$ consists of $s$ singular points of
length 3, thus by Lemma \ref{2.4}(1),
$$\RHS\le \frac{1}{4}+0-\frac{4}{13} < 0,$$
a contradiction.

\medskip\noindent
Case 10. $c_2(S')=3s+4$ and $K_{S'}^2=8-3s$ ($s\ge 3$).

\medskip\noindent
In this case  $|\Sing(S)|\ge (b_2(S')-1)/3=(3s+1)/3$, hence
$|\Sing(S)|\ge s+1$.

Since $|\Sing(S)|\le 4$, $s=3$ and $|\Sing(S)|=4$.
Thus
$$\LHS= (b_2(S')-1)+|\Sing(S)|-c_2(S')+\frac{1}{3}K_{S'}^2
= \frac{2}{3}.$$
Since $b_2(S')-1=1+3+3+3=2+2+3+3$, $\Sing(S)$ consists either of 1 singular point of
length 1 and $3$ singular points of
length 3, or 2 singular points of
length 2 and $2$ singular points of
length 3, thus by Lemma \ref{2.4}(1),
$$\RHS\le \frac{1}{2}+0\quad {\rm or}\quad \le \frac{2}{9}+
\frac{1}{4}+0
\quad {\rm or}\quad \le\frac{1}{3}+ \frac{1}{15}+\frac{1}{4}+0,$$
all smaller than $\frac{2}{3}$,
a contradiction.

\medskip\noindent
Case 11. $c_2(S')=3s+2$ and $K_{S'}^2=10-3s$ ($s\ge 4$).

\medskip\noindent
In this case  $|\Sing(S)|\ge (b_2(S')-1)/3=(3s-1)/3$, hence
$|\Sing(S)|\ge s$.

Since $|\Sing(S)|\le 4$, $s=4$ and $|\Sing(S)|=4$.
Thus
$$\LHS= (b_2(S')-1)+|\Sing(S)|-c_2(S')+\frac{1}{3}K_{S'}^2
= \frac{1}{3}.$$
Since $b_2(S')-1=2+3+3+3$, $\Sing(S)$ consists of 1 singular point of
length 2 and $3$ singular points of
length 3, thus by Lemma \ref{2.4}(1),
$$\RHS\le \frac{1}{3}+ \frac{1}{4}+0-\frac{4}{13} < \frac{1}{3},$$
a contradiction.

This completes the proof of Theorem \ref{main}.

\begin{corollary}\label{cor} Koll\'ar's conjecture holds true
if in addition the exceptional divisor in
a minimal resolution of $S$ has at most 3 components over each
singular point of $S$.
\end{corollary}

\section{Examples and Further Discussion}

\begin{example}\label{Ishida} In \cite{Ishida} M. Ishida discusses an elliptic surface $Y$ with
$p_g=q=0$ with two multiple fibres, one of multiplicity 2 and one
of muliplicity 3, and proves that the Mumford fake plane is its
cover of degree 21, non-Galois. The surface $Y$ is a Dolgachev
surface \cite{BHPV}. In particular, it is simply connected and of
Kodaira dimension 1. Besides the two multiple fibres, its elliptic
fibration $|F_Y|$ has 4 more singular fibres $F_1, F_2, F_3, F_4$,
all of type $I_3$. It has also a sixtuple section $E$ which is a
$(-3)$-curve meeting one component of each of $F_1, F_2, F_3$  in
6 points, and two components of $F_4$ in 1 point and 5 points
each. One can contract 9 rational curves to get a singular surface
$S$ with 3 singular points of type $\frac{1}{3}(1,2)$ and one
singular point of type $\frac{1}{7}(1,3)$. The surface $S$
satisfies the condition (1) and (3) of Theorem \ref{main}, but not
(2). Indeed, $H_1(S^0,\bbZ)=\bbZ/3\bbZ$.
\end{example}

\begin{example}\label{ex} It is shown \cite{Ke} that there is another Dolgachev
surface $X$ which is birational to a cyclic cover of degree 3 of
Ishida surface $Y$. On $X$ there are 9 smooth rational curves
forming a configuration
$$(-2)\textrm{---}(-2)\textrm{---}(-3) \quad (-2)\textrm{---}(-2)\textrm{---}(-3) \quad (-2)\textrm{---}(-2)\textrm{---}(-3)$$
which can be contracted to $3$ singular points of type
$\frac{1}{7}(1,3)$. The resulting singular surface $S$ satisfies
the condition (1) and (3) of Theorem \ref{main}, but not (2). In
this case, $H_1(S^0,\bbZ)=\bbZ/7\bbZ$.
\end{example}

Finally we consider surfaces $S$ with rational double points only.

\begin{proposition}\label{rdp} Let $S$ be a singular projective surface with
rational double points such that $H^2(S, \bbQ)\cong \bbQ$ and
$H_1(S^0,\bbZ)=0$. Let $f: S'\to S$ be a minimal
resolution. Then one of the following cases occurs.
\begin{itemize}
\item[(1)] $S$ is rational.
\item[(2)] $S'$ is a minimal surface of general type,
with $q=p_g=0$,  $H_1(S',\bbZ)=0$, and
\begin{itemize}
\item[(2-1)] $K_{S'}^2=1$ and $\calR\cong E_8$, or
\item[(2-2)] $K_{S'}^2=2$ and $\calR\cong E_7$, or
\item[(2-3)] $K_{S'}^2=3$ and $\calR\cong E_6$, or
\item[(2-4)] $K_{S'}^2=4$ and $\calR\cong D_5$, or
\item[(2-5)] $K_{S'}^2=5$ and $\calR\cong A_4$.
\end{itemize}
\end{itemize}
\end{proposition}

\begin{proof} Since $S$ has
rational double points only, $f^*K_S=K_{S'}$. If $K_S$ is
anti-numerically effective, so is $K_{S'}$, hence $S'$ is rational.
If $K_S$ is numerically effective, so is $K_{S'}$, hence $S'$ is
minimal.
We use Proposition \ref{3}. We need to rule out the second possibility from
Proposition \ref{3}. Suppose that the second case occur. Since $K_{S'}$ is a
rational multiple of a fibre of the elliptic fibration, the
exceptional divisor of $f$ is supported in a union of fibres. This
contradicts to $\kappa(S^0)=2$.

Next Assume that $S'$ is a minimal surface of general type.
The divisor  $K_{S'}$ is an integral divisor belonging to
$\calR^{\perp}$, hence $K_{S'}=mv$ for some integer $m$,
where $v$ is a generator of $\calR^{\perp}$.
Since $v^2=|\det\calR|$, $K_{S'}^2=m^2|\det\calR|$.
This leaves the five cases (2-1)-(2-5) and two more
\begin{itemize}
\item[(2-6)] $K_{S'}^2=6$ and $\calR\cong A_1\oplus A_2$, or
\item[(2-7)] $K_{S'}^2=8$ and $\calR\cong A_1$.
\end{itemize}
Both are ruled out by Theorem \ref{main}.
\end{proof}

\begin{remark} If one loosens the bound to 4 on the number of components in
Condition (3) of Theorem \ref{main}, one already encounters a non-trivial
problem to rule out the possibility (2-5) from the above proposition.
\end{remark}


\end{document}